\documentclass{smfart}

\usepackage[latin1]{inputenc} \usepackage[T1]{fontenc} \usepackage{lmodern}

\usepackage{amssymb,amsmath,latexsym,amsthm}

\newtheorem{lemma}{Lemma}[section]
\newtheorem{proposition}{Proposition}[section]
\newtheorem{theorem}{Theorem}

\newtheorem*{definition}{Definition}

\numberwithin{equation}{section}

\usepackage{algorithm}
\usepackage{algorithmic}

\newcommand{\ALRETURN}{\item[\ \ \ \textbf{return}]}

\input xy
\xyoption{all}

\def\Ot{\mathop{{\tilde{O}}}\nolimits }
\def\bK{{\bf K  }}
\def\bL{{\bf L  }}

\let\set\mathbb
\def\FF{{\set F}}
\def\PP{{\set P}}
\def\ZZ{{\set Z}}
\def\QQ{{\set Q}}
\def\NN{{\set N}}

\def\pod{\mathop{\rm{PDiv}}\nolimits }
\def\loss{\mathop{\rm{Loss}}\nolimits }
\def\lploss{\mathop{\rm{LpLoss}}\nolimits }

\begin{document}

\title[Elkies subgroups of elliptic curve $\ell$-torsion points]{On Elkies
  subgroups of $\ell$-torsion points in elliptic curves defined over a finite
  field}

\author{{\sc Reynald} LERCIER}
\address{
  \textsc{DGA}/\textsc{C\'ELAR}, La Roche
  Marguerite, F-35174 Bruz, and\\
  IRMAR, Universit\'e de Rennes 1, Campus
  de Beaulieu, F-35042 Rennes}
\email{reynald.lercier@m4x.org}
\urladdr{http://perso.univ-rennes1.fr/reynald.lercier/}

\author{{\sc Thomas} SIRVENT}
\address{
  \textsc{DGA}/\textsc{C\'ELAR}, La Roche
  Marguerite, F-35174 Bruz, and\\
  IRMAR, Universit\'e de Rennes 1, Campus
  de Beaulieu, F-35042 Rennes}
\email{thomas.sirvent@m4x.org}
\urladdr{http://perso.univ-rennes1.fr/thomas.sirvent/}

\begin{abstract}
  En sous-résultat de l'algorithme de Schoof-Elkies-Atkin pour compter le
  nombre de points d'une courbe elliptique définie sur un corps fini de
  caractéristique $p$, il existe un algorithme qui, pour $\ell$ 
  un nombre premier d'Elkies, calcule des points de $\ell$-torsion dans une
  extension de degré $\ell-1$ à l'aide de $\Ot(\ell \, \max(\ell,
  \log q)^2)$ opérations élémentaires à condition que
  $\ell\leqslant p/2$.

  Nous combinons ici un algorithme rapide dû à Bostan, Morain, Salvy et
  Schost avec l'approche $p$-adique suivie par Joux et Lercier pour
  obtenir pour la première fois un algorithme valide sans limitation
  sur $\ell$ et $p$ et de complexité similaire.
\end{abstract}

\begin{altabstract}
  As a subproduct of the Schoof-Elkies-Atkin algorithm to count points on
  elliptic curves defined over finite fields of characteristic $p$, there
  exists an algorithm that computes, for $\ell$ an Elkies prime,
  $\ell$-torsion points in an extension of degree $\ell-1$ at cost $\Ot(\ell
  \, \max(\ell, \log q)^2)$ bit operations in the favorable case where
  $\ell\leqslant p/2$.

  We combine in this work a fast algorithm for computing isogenies due to
  Bostan, Morain, Salvy and Schost with the $p$-adic approach followed by Joux
  and Lercier to get for the first time an algorithm valid without any
  limitation on $\ell$ and $p$ but of similar complexity.
\end{altabstract}

\maketitle
\let\languagename\relax

\smallskip

\section{Introduction}

Let $\bK$ be a finite field with $q$ elements and $E$ be an elliptic
curve over $\bK$ given by a plane equation of the form
\begin{equation}\label{eq:wec}
  y^2+a_1xy+a_3y = x^3+a_2x^2+a_4x+a_6
\end{equation}
where the coefficients $a_1$, $a_2$, $a_3$, $a_4$ and $a_6$ are elements of
$\bK$. For any field $\bL$ such that $\bK\subset \bL$, we denote by $E(\bL)$
the set of $\bL$-points of $E$, \textit{i.e.} the set of solutions in $\bL$ of
Equation~(\ref{eq:wec}), plus the additional point at infinity $O$ with
homogeneous coordinates $(0:1:0)$. The curve $E/\bK$ has a structure of
commutative algebraic group with neutral element $O$, derived from the secant
and tangent rules. Its order is equal to $q+1-t$ with $t \in \ZZ$ such that
$|t|\leqslant 2\sqrt{q}$. \smallskip

We are interested in the determination of $\ell$-torsion points of $E$, that
is the set $E[\ell]$ of points $P$ of $E(\overline{\bK})$ such that $\ell P =
O$ for prime integers $\ell$, distinct from $p$.  This group is isomorphic to
$\ZZ / \ell \ZZ \times \ZZ / \ell \ZZ$ (cf.\cite[p.  89]{00079}), its cardinal
is thus $\ell^2$. In fact, the multiplication by $\ell$ is given by a rational
transformation of $\PP^2(\bK)$, of degree $\ell^2$, of the form $(x:y:z)
\mapsto (X_\ell(x,y,z) : Y_\ell(x,y,z) : Z_\ell(x, y, z))$ where $X_\ell$,
$Y_\ell$ and $Z_\ell$ are three homogeneous polynomials of degree $\ell^2$ and
$\ell$-torsion points are explicitely given by $Z_\ell(x, y, z) = 0$. Excluding $O$, this
equation can be easily transformed into an equality of the form $f_\ell(x) = 0$
where $f_\ell$ is an univariate polynomial of degree $(\ell^2-1)/2$, called
the $\ell$-th division polynomial. \smallskip

The improvements by Atkin and Elkies to Schoof's algorithm for counting points
on elliptic curve stem from the fact that when the principal ideal $(\ell)$
splits in the imaginary quadratic field $\QQ(\sqrt{t^2-4q})$, in half
the cases thus, there exists two subgroups of degree $\ell$ in $E[\ell]$
defined in a degree $\ell-1$ extension of $\bK$. Such an integer $\ell$ is
called an Elkies prime.  In this work, we more precisely focus on algorithmic
efficient ways to compute degree $(\ell-1)/2$ polynomials the roots of which
are abscissas of points contained in such subgroups. We call these subgroups,
and these degree $(\ell-1)/2$ polynomials over $\bK$,
$\ell$-th Elkies subgroups, and
$\ell$-th Elkies polynomials.\medskip

Our main result, where we classically denote by $\phi_1(x_1$, $\ldots, x_k) =
\Ot(\phi_2(x_1,$ $\ldots, x_k))$ functions $\phi_1$ and $\phi_2$ such that there
exists an integer $k$ with $\phi_1(x_1, \ldots, $ $x_k) = O ( \phi_2(x_1, \ldots,
x_k) \, \log^k  \phi_2(x_1, \ldots, x_k)  )$, is as follows.
\begin{theorem}
  Let $E$ be an elliptic curve defined over a finite field $\bK$ with $q$
  elements and $\ell$ be an Elkies prime, distinct from the characteristic
  of $\bK$, then there exists an algorithm which computes an $\ell$-th Elkies
  polynomial at cost $\Ot(\ell \, \max(\ell, \log q)^2)$ bit operations and
  space.
\end{theorem}

This problem is closely related to the problem of computing separable
isogenies of degree $\ell$ between two elliptic curves since an application of
Velu's formulas~\cite{00103} with starting point such polynomials yields an
isogeny.  Especially, counting points on elliptic curves first raised interest
for such computations. But isogenies now play a role in numerous other fields,
for instance to protect elliptic curve cryptographic devices against physical
side attacks~\cite{Smart2003}, to improve Weil descent to calculate elliptic
discrete logarithms~\cite{GaHeSm2002}, to decrease the complexity of computing
discrete logarithms in some family of finite fields~\cite{CouLer2007}, to
exhibit normal basis in finite field extensions~\cite{CouLer2008}, etc.

\medskip  

We first recall in Section~\ref{sec:related} the complexity of the algorithms
known to solve this problem. In Section~\ref{sec:largechar}, we focus on the
fastest algorithm in finite fields of large characteristic published so far,
due to Bostan, Morain, Salvy and Schost~\cite{00087}. We then show in
Section~\ref{sec:smallchar} how we can combine this algorithm with the
$p$-adic approach introduced by Joux and Lercier in~\cite{00088} to get a fast
algorithm in any finite field and we clarify that we need a $p$-adic precision
of only $O(\log^2 \ell/\log p)$.  A detailed example is
given in Section~\ref{sec:experiments}.

\section{Related work}
\label{sec:related}

For the sake of simplicity, we restrict ourself to finite fields $\bK$ of
characteristic larger than three, and to prime integers $\ell > 2$.
In this case, an elliptic curve is
simply given by a plane equation of the form $y^2 = x^3 + a_4x +a_6$. Its
discriminant, always non zero, is equal to $\Delta_E = -16 (4 {a_4}^3 + 27
{a_6}^2)$ and its $j$-invariant is equal to $j_E = - 12^3 \, {(4
  a_4)^3}/{\Delta_E}.$

\subsection{Naive approach}
$\ell$-th Elkies polynomials are factors of the $\ell$-th division polynomial
$f_\ell$. Therefore, a naive approach consists in computing $f_\ell$, which
can be done at cost $\Ot (\ell^2 \, \log q)$ elementary operations thanks to
a ``Square and Multiply'' method~\cite{00079}, and then in factorizing it with
cost $\Ot (\ell^{1,815 \times 2} \, \log^2 q)$~\cite{00105}. This algorithm
needs a total of $\Ot (\ell^4 \, \log^2 q)$ bit operations.
    
\subsection{Schoof-Elkies-Atkin framework}

Let $\pi_E$ be the Frobenius endomorphism of $E$. Its restriction to
$E[\ell]$, seen as a $\FF_\ell$-vector space of dimension two, is still an
endomorphism. When $\ell$ is an Elkies prime, its eigenspaces correspond to
$\ell$-th Elkies subgroups $C$ of $E[\ell]$ and from each $C$ one can construct
an isogeny of degree $\ell$ between $E$ and the elliptic curve $E'=E/C$,
defined over $\bK$. \medskip

The following algorithm takes advantage of these facts.
\begin{description}
\item[Step 1] Compute the modular polynomial of degree $\ell$, $\Phi_\ell(X, Y)$,
  equation of the modular curve $X_0(\ell)$. This is a bivariate symmetric
  polynomial, of degree $\ell + 1$ in $X$ and $Y$, whose coefficients are
  integers of $\Ot(\ell)$ bits (\textit{cf.}~\cite{Cohen84}). $j$-invariants of
  $\ell$-isogenous elliptic curves are roots of $\Phi_\ell(X, Y)$.
\item[Step 2] Compute roots $j_1$ and $j_2$ of $\Phi_\ell(X, j_E)$.
\item[Step 3] Compute a normalized Weierstrass equation for elliptic curves of
  $j$-invariants $j_1$ and $j_2$, and the sum $p_1$ of the abscissas of points
  in the kernel of the isogeny, using the polynomials $\Phi_\ell$, $\partial
  \Phi_\ell / \partial X$, $\partial \Phi_\ell /
  \partial Y$, $\partial^2 \Phi_\ell / \partial X^2$, $\partial^2 \Phi_\ell /
  \partial X \partial Y$, $\partial^2 \Phi_\ell / \partial Y^2$
  (\textit{cf.}~\cite{00107}).
\item[Step 4] Compute from each isogenous curve, a $\ell$-th Elkies polynomial
  thanks to the kernel of the corresponding isogeny.
\end{description}
\medskip

The complexity of the method comes now.
\begin{description}
\item[Step 1] The modular polynomial $\Phi_\ell(X, Y)$ has $O(\ell^2)$
  coefficients, each with about $\Ot(\ell)$ bits. There exists methods to
  compute this polynomial at cost quasi-linear in its size, \textit{i.e.}
  $\Ot(\ell^3)$ bit operations (\textit{cf.}~\cite{00104}). We need to reduce this polynomial modulo $p$,
  that is $\Ot(\ell^3)$ bit operations too. The result is then of size
  $O(\ell^2 \, \log p )$ bits.
\item[Step 2] With the help of Horner's method, the evaluation of
  $\Phi_\ell(X, Y)$ at $j_E$ costs $\Ot(\ell^2 \, \log q)$ bit operations. In
  order to compute roots of the resulting degree $\ell+1$ polynomial, we have
  first to compute its gcd with $X^q - X$, that is $\Ot(\ell \, \log^2 q)$
  bit operations (\textit{cf.}~\cite{00123}). We obtain a degree $2$
  polynomial whose roots can then be found with negligible cost.
\item[Step 3] The computations of the derivatives of $\Phi_\ell$ and their
  evaluations can be done at cost $\Ot(\ell^2 \, \log q)$ bit operations.
\item[Step 4] Here, we have to distinguish several cases.
  \begin{itemize}
  \item In finite fields of large characteristic, the best algorithm known so
    far to compute isogenies is due to Bostan et al.~\cite{00087} and takes
    time $\Ot(\ell \, \log q)$ bit operations.
  \item In finite fields of small but \textit{fixed} characteristic, the best
    algorithm known is due to Couveignes~\cite{00114} and needs $\Ot(\ell^2
    \log q)$ bit operations (but the contribution of $p$ in the $\Ot$
    complexity constant is exponential in $\log p$).
  \item In between, that is finite fields of small but \textit{non-fixed}
    characteristic, the best algorithm is due to Joux and Lercier~\cite{00088}
    and needs $\Ot((1 + \ell/p) \, \ell^2 \log q)$ bit operations.
  \end{itemize}
\end{description}
\medskip

The best total complexity is thus equal to $\Ot(\ell \, \max(\ell,
\log q)^2)$, achieved in finite fields of large characteristic. But, in
finite fields of small characteristic, the complexity can be as large as
$\Ot(\ell^3\log q)$ bit operations when $\ell \gg p$. \medskip

This work yields an algorithm of same complexity as in the large
characteristic case without any limitation on the characteristic or the degree
of the base field $\bK$.

\section{The large characteristic case}
\label{sec:largechar}
  
In order to get an algorithm with good complexity in finite fields of small
characteristic too, we first reformulate the algorithm of Bostan et al. in
such a way that its extension in the $p$-adics is more easily studiable.  The
general strategy is the same except that we take into account some
specificities of the involved differential equation in the resolution. As a
result, we obtain a precise and compact algorithm (cf.
Algorithm~\ref{algo:eqdiff}).

\subsection{Differential equation}

In a field $\bK$ of characteristic larger than three, an isogeny between two
elliptic curves, $E:y^2 = x^3 + a_4 \, x + a_6$ and $E':y^2 = x^3 + a'_4 \, x
+ a'_6$, can be given by
\begin{displaymath}
  I(x, y) = \left( \frac{N(x)}{D(x)}, c y \left( \frac{N(x)}{D(x)} \right)' \right)\,,
\end{displaymath}
where $N$ and $D$ are unitary polynomials of degree $\ell$ and $\ell-1$.  When
$c$ is equal to one, the isogeny is said to be normalized. This is in particular the case in the Schoof-Elkies-Atkin framework.

If we now state that the image of a point of $E$ by $I$ is on $E'$, we get the
following differential equation
\begin{equation}
  \label{eq:isogenie_eqdiff}
  (x^3 + a_4 \, x + a_6) \left( \frac{N(x)}{D(x)} \right)'^{\, 2} = \left( \frac{N(x)}{D(x)} \right)^3 + a'_4 \left( \frac{N(x)}{D(x)} \right) + a'_6.
\end{equation}
This equation can be solved with a Taylor series expansion of $N(x)/D(x)- x$ in $1/x$ for $1/x$ close to $0$. The relations obtained thanks to
Equation~(\ref{eq:isogenie_eqdiff}) enable to compute by recurrence each
coefficient in turn, if the first coefficients are known. It is then possible
to recover $N$ and $D$ with the help of Berlekamp-Massey's algorithm, or one
of its optimized variant. In~\cite{00087}, one takes advantage of a Newton
algorithm so that the number of coefficients computed at each iteration
doubles.

More precisely, let $S$ be defined by
\begin{displaymath}
  S(x) = \sqrt{\frac{D \left( 1 / x^2 \right)}{N \left( 1 / x^2 \right)}} \, ,
  \, \text{ or equivalently } \, \frac{N(x)}{D(x)} = \frac{1}{S \left( {1}/{\sqrt{x}} \right)^2}.
\end{displaymath}
Equation~(\ref{eq:isogenie_eqdiff}) becomes
\begin{displaymath}\label{eq:eqdiff}
  (a_6 \, x^6 + a_4 \, x^4 + 1) S' (x)^2 = 1 + a'_4 \, S(x)^4 + a'_6 \, S(x)^6.
\end{displaymath}
At the infinity, $N(x)/D(x)$ has a series expansion of the form $x +
O(1)$. We thus have $S(x) = x + O(x^3)$ and this knowledge is finally enough to
completely recover $N(x)/D(x)$.
    
\subsection{Resolution}
    
We consider more generally equations of the form $S'^2 = G \cdot (H \circ S)$.
In Equation~(\ref{eq:eqdiff}), we have for instance $H(z) = a'_6 \, z^6 + a'_4
\, z^4 + 1$ and $G(x) = 1/(a_6 \, x^6 + a_4 \, x^4 + 1)$.  We now look for a
solution modulo $x^\mu$, where $\mu$ is any integer given in input.  The way
to solve this equation is first to assume that we know the solution modulo
$x^d$ and then, thanks to a Newton iteration, to obtain a solution modulo
$x^{2d}$.  After roughly $\log \mu$ such iterations, one gets the full
solution.

We now present a compact algorithm for this task. Its complexity can be easily
established, it is equal to $\Ot(\mu \log q)$ bit operations. Its correctness
is slightly more difficult to prove and we delay it to Appendix~\ref{sec:proof}.
\begin{algorithm}
  \begin{algorithmic}
    \REQUIRE $\mu \in \NN$, $(\alpha, \beta) \in \bK^2$, $H \in \bK[z]$, $G \in \bK[[x]]$
    \ENSURE $S \in \bK[x]$, a solution of the differential equation modulo $x^\mu$
    \smallskip
    \STATE $d \longleftarrow 2$, \, $U \longleftarrow 1 / \beta$, \, $J \longleftarrow 1$, \, $V \longleftarrow 1$
    \STATE $S \longleftarrow \alpha + \beta \, x + \big[ \big( G'(0) + H'(\alpha) \, \beta^3 \big) / (4 \beta) \big] \, x^2$
    \WHILE {$(d < \mu - 1)$}
    \STATE $U \longleftarrow U \cdot \big( 2 - S' \cdot U \big) \bmod x^d$
    \STATE $V \longleftarrow \big( V + J \cdot (H \circ S) \cdot ( 2 - V \cdot J) \big) \, / \, 2 \bmod x^d$
    \STATE $J \longleftarrow J \cdot \big( 2 - V \cdot J \big) \bmod x^d$
    \STATE $S \longleftarrow S + V \cdot \int \big( G \cdot (H \circ S) - S'^2 \big) \, \big( U \cdot J \, / \, 2 \big) \, dx \bmod x^{\min(2d+1,\mu)}$
    \STATE $d \longleftarrow 2d$
    \ENDWHILE
    \ALRETURN $S$
  \end{algorithmic}
  \caption{Solving equation $S'^2 = G \cdot (H \circ S)$, $S(0) = \alpha$ and $S'(0) = \beta$.}
  \label{algo:eqdiff}
\end{algorithm}
    
\begin{proposition}\label{prop:eqdiff}
  Let $(\alpha, \beta) \in \bK^2$ where $\bK$ is a finite field of
  characteristic $p$, let $G$ be a formal series defined over $\bK$, let $H$ be
  a polynomial defined over $\bK$ such that $H(\alpha) = 1$ and $\beta^2 =
  G(0) \neq 0$. Let $\mu \in \{1, \ldots, p\}$, then Algorithm~\ref{algo:eqdiff} computes
  a Taylor series (modulo $x^\mu$) of the solution $S$ of the differential
  equation
  \begin{displaymath}
    S'(x)^2 = G(x) \, H(S(x)) \, , \, S(0) = \alpha \, , \, S'(0) = \beta.
  \end{displaymath}
\end{proposition}
    
    \smallskip
    
\subsection{Full algorithm}
\label{isogenie_calcul}
    
We first compute $G(x) = 1/(a_6 \, x^6 + a_4 \, x^4 + 1)$ modulo $x^{4 \ell -
  1}$ thanks to the classical iterative following formula, $G_1(x) = 1$,
$G_{2d}(x) = G_d(x) \left(2 - G_d(x) \cdot (a_6 \, x^6 + a_4 \, x^4 + 1)
\right) \bmod x^{2d}$.  We then apply Algorithm~\ref{algo:eqdiff} to $G(x)$
and $H(z) = a'_6 \, z^6 + a'_4 \, z^4 + 1$ with $\mu = 4 \ell$, $\alpha = 0$ and $\beta = 1$.
    
The obtained solution $S$ is odd, we define from it
\begin{displaymath}
  T(x) = \sum_{i=0}^{2 \ell - 1} t_i \, x^i \, , \text{ where } \forall i \in \{0, \ldots, 2 \ell - 1\},\ t_i = s_{2i+1}.
\end{displaymath}
    
We denote by $R(x)$ the inverse of the square of $T(x)$, modulo $x^{2 \ell}$,
with the same inverse formulas as those used for $G$. We then have
\begin{displaymath}
  \frac{N(x)}{D(x)} = x \, R \left( \frac{1}{x} \right)\,, \text{ \textit{i.e.} } R(x) = \frac{x^\ell \, N(1/x)}{x^{\ell-1} \, D(1/x)}\,.
\end{displaymath}

Applying Berlekamp-Massey algorithm~\cite{00119,00118,00120} or one of its
optimized variant~\cite{00122,00121} to $R$ yields $D$ and the searched
$\ell$-th Elkies polynomial is equal to the square root of $D$.
    
\section{Extension to any finite field}
\label{sec:smallchar}

To extend the Schoof-Elkies-Atkin framework in any characteristic, the
techniques developed in~\cite{00088} give the general idea: to use the
$p$-adics to authorize divisions by the characteristic $p$ of the field. These
divisions make it possible to use in any finite field algorithms primarily
designed in large characteristic. There exists one main obstacle with this
approach. Calculations in the $p$-adics imply losses of precision at the time
of divisions by $p$. It is thus necessary to anticipate a sufficient precision,
which results in an increase in the size of the handled objects.

One could hope to perform this lift in the $p$-adics only in the last stage of
the algorithm, \textit{i.e.} for the calculation of the isogeny.  It is
actually not possible because fast algorithms for computing isogenies need
normalized models for the isogenous curves.

It is thus necessary to lift in the $p$-adics from the very beginning of the
algorithm. It is exactly what is done in~\cite{00088}, with a $p$-adic
precision linear in $\ell$.  Instead, we consider here the techniques of
\cite{00087}, and one shows that the necessary $p$-adic precision can be
brought back to only $O(\log^2 \ell/\log p)$. The total complexity of the
algorithm is then similar to the one of the large characteristic case, that is
$\Ot(\ell \, \max(\ell, \log q)^2)$.

\subsection{Lifting curves and isogenies}

One starts by lifting arbitrarily the curve $E$ in the $p$-adics. Any
coefficient $\bar{a}_4$ and $\bar{a}_6$ such that $\bar{a}_4 = a_4 \bmod p$ and
$\bar{a}_6 = a_6 \bmod p$ is appropriate and one works on the elliptic
curve $\bar{E}/\QQ_{q}$ with model $y^2 = x^3 + \bar{a}_4 \, x +
\bar{a}_6$.

The computation of the $j$-invariant $\bar{j_E}$ of the curve $\bar{E}$, of
the solutions $\bar{j_1}$ and $\bar{j_2}$ of the equation $\phi_\ell(x,
\bar{j_E}) = 0$, as well as Weierstrass models of the corresponding curves
$\bar{E}_1$ and $\bar{E}_2$, proceeds exactly as in the SEA framework.
The curves $\bar{E}_1$ and $\bar{E}_2$ are $\ell$-isogenous with the curve
$\bar{E}$, and the isogenies can be calculated as in the large characteristic
case. 

Projection $E_1$ of the curve $\bar{E}_1$ on the base field $\bK$
is $\ell$-isogenous with $E$, and the connecting isogeny is the projection on
the base field of the isogeny connecting $\bar{E}$ to $\bar{E}_1$. It is the
same for $E_2$. It is thus enough to project the denominators of the
isogenies on $\bK$ to identify the required factors of the
$\ell$-th division polynomial of $E$.

\subsection{$p$-adic computations}

From now on, we are interested in the $p$-adic precision of the lift of the
elliptic curve $E$. This precision must be large enough so that at the end of
the resolution of the differential equation with Algorithm~\ref{algo:eqdiff},
the result $S$ can be reduced in $\bK$.\smallskip

To this purpose, we need first some definitions.
\begin{definition}
  For any positive integer $r$, one defines $\pod(p, r)$ by the
  largest power of $p$ which divides $r$,
  \begin{math}
    \pod(p, r) = \max \big\{ k \in \NN \, | \, p^k \,
    \text{ divides } \, r \big\}\,.   
  \end{math}

  We denote by $\loss(p, \ell)$ the sum
  \begin{math}
    \sum_{1 \, \leqslant \, i \, < \, \log_2(4 \ell - 1)} \lploss (p, \ell, i),
  \end{math}
  where
  \begin{multline*}
    \lploss(p, \ell, i) \, = \max \, \big\{
    \pod(p, r) \, |
    \, 2^i + 1 \leqslant r \leqslant \min(2^{i+1}, 4 \ell - 1) \big\}.
  \end{multline*}
\end{definition}
\medskip

The following lemma relates the precision needed to the function $\loss$.
\begin{lemma}
  \label{isogenie_theo}
  Let $\mu$ be the $p$-adic precision of the coefficients $\bar{a}_4$ and
  $\bar{a}_6$, then when $\mu > \loss(p, \ell)$ the polynomials  $U$, $V$,
  $J$ and $S$ computed in Algorithm~\ref{algo:eqdiff} have $p$-adic integer
  coefficients. Furthermore the precision of the result $S$ is at least equal
  to $(\mu - \loss(p, \ell))$. 
\end{lemma}

\begin{proof}
  One proves this theorem by recurrence on $j$, the number of iterations of
  the loop ``while'' in Algorithm~\ref{algo:eqdiff}. We assume that at rank
  $j$, $0 \leqslant j < \log_2(4 \ell - 1)$, the polynomials $U$, $V$, $J$ and
  $S$ have $p$-adic integer coefficients and that their precision is at
  least equal to $\mu - \sum_{1 \, \leqslant \, i \, \leqslant \, j} \lploss
  (p, \ell, i)$.

  \noindent
  \textsc{Initialization.}
  In input of the algorithm, we have $\alpha = 0$, $\beta = 1$, $H(z) =
  \bar{a}'_6 \, z^6 + \bar{a}'_4 \, z^4 + 1$ and $G(x) = 1/(\bar{a}_6 \, x^6 +
  \bar{a}_4 \, x^4 + 1)$. The elements $\bar{a}_4$, $\bar{a}_6$, $\bar{a}'_4$
  and $\bar{a}'_6$ are integers of precision $\mu$ and thus $G$ and $H$ are
  of precision $\mu$ too (no division by $p$ occurs in the computation of
  $G$).  The same is true for $U$, $V$, $J$ and $S$.  
   
  \noindent
  \textsc{Heredity.}
  Let $j < \log_2(4 \ell - 1)$, we suppose the assumption true at rank $j-1$.
  At the $j^\text{th}$ iteration, polynomials $U$, $V$ and $J$ are updated via
  additions, multiplications, derivations and compositions of the values of
  $U$, $V$, $J$ and $S$ before the entry in the loop. All these operations
  preserve the precision, polynomials $U$, $V$ and $J$ have thus $p$-adic
  integer coefficients with precision at least equal to $\mu - \sum_{1 \,
    \leqslant \, i \, \leqslant \, j-1} \lploss (p, \ell, i)$.

  For $S$, except the integral operation, the calculations preserve the
  precision.  Coefficients of the series after the integral operation are
  inverses of degrees between $2^j + 1$ and
  $\min(2^{j+1}, 4 \ell - 1)$.  The largest power of $p$ by which we carry out
  a division is thus $\lploss(p, \ell, j)$. The absolute precision of the
  coefficients of $S$ is thus higher or equal to $\mu - \sum_{1 \, \leqslant
    \, i \, \leqslant \, j} \lploss (p, \ell, i)$. Furthermore, since this precision is positive, each coefficient of $S$ is a lift of the coefficient of the series deduced from the isogeny over $\bK$, and these coefficients are $p$-adic integers.
\end{proof}
    
To minimize the loss of precision, we may use the additional fact that $S$ is
odd. We thus have to consider only coefficients of odd degree in the
algorithm and the loss of precision in the loop of the algorithm becomes
\begin{displaymath}
  \lploss \, '(p, \ell, i) \, = \max \, \big\{ \pod(p, 2 r + 1) \, \big/ \,
  2^{i-1} \leqslant r \leqslant \min(2^{i} - 1, 2 \ell - 1) \big\}\,.
\end{displaymath}
    
Lemma~\ref{isogenie_prop} yields a clear asymptotic bound on the loss of
precision stated in Lemma~\ref{isogenie_theo}.
\begin{lemma}
  \label{isogenie_prop}
  We have $\loss(p, \ell) = O \left( \log^2 \ell/\log p \right)$.
\end{lemma}
\begin{proof}
  For all $i < \log_2(4 \ell - 1)$, $\lploss(p, \ell, i)$ is the largest power
  of $p$ which divides a range of integers, at most equal to $2^{i+1}$, we
  have therefore $\lploss(p, \ell, i) \leqslant \log_p 2^{i+1}$, and
  \begin{equation*}
    \begin{array}{rcl}
      \loss(p, \ell) & \leqslant & \log_p 2 \, \left( \sum_{1 \, \leqslant \, i \, < \,
          \log_2(4 \ell - 1)} \, (i + 1) \right)\,, \\
      & \leqslant & \log_p 2 \,\log_2(4 \ell - 1)\left(
        \log_2(4 \ell - 1) + 1 \right) \,, \\
      %
      %
      & \leqslant & \left( \log_2(4 \ell - 1) + 1 \right)^2/\log_2 p\,. \\
    \end{array}
  \end{equation*}      
\end{proof}
    
We finally can state our main result.
\begin{proposition}\label{prop:mainresult}
  A $p$-adic precision of $O(\log^2 \ell/\log p)$ is asymptotically enough
  to compute a $\ell$-th Elkies polynomial. The total computation needs
  $\Ot(\ell \, \max(\ell, \log q)^2)$ bit operations.
\end{proposition}

\begin{proof}
  Computations performed in the Schoof-Elkies-Atkin framework, especially
  calls to Algorithm~\ref{algo:eqdiff}, are realized in the $p$-adics with
  precision at most $O(\log^2 \ell/\log p)$. This precision does not modify
  the $\Ot$ complexities of the large characteristic case and we still have in
  the $p$-adic case a total complexity equal to $\Ot(\ell \, \max(\ell,$
  $\log q)^2)$ bit operations.
\end{proof}
    
\section{Experiments}
\label{sec:experiments}

We have implemented this algorithm in the computer algebra system {\sc
  magma}.  Thanks to it, we were able to observe that the bound on the
precision stated in Proposition~\ref{prop:mainresult} is tight. We can
illustrate the method with an example too.

\subsection{$p$-adic precision}

Figure~\ref{fig:isogenie_precision} shows the evolution of the precision when
$p$ and $\ell$ vary. The ``The(oretical)'' bound mentioned corresponds to
$\loss(p, \ell)$ calculations. The ``Obs(erved)'' bound is what seems
necessary at the time of calculations (checked on some examples).
     
\begin{figure}[htbp]
  \centering
  \begin{tabular}{ccc}
    \begin{tabular}{|c|c|c|}
      \hline
      \multicolumn{3}{|c|}{$p=5$} \\
      \hline
      $\ell$ & Obs. & The. \\\hline\hline
      7 & 5 & 5 \\\hline
      11 & 6 & 6 \\\hline
      13 & 6 & 7 \\\hline
      17 & 7 & 8 \\\hline
      19-31 & 8 & 9 \\\hline
      37 & 11 & 11 \\\hline
      41-61 & 11 & 12 \\\hline
      67 & 13 & 14 \\\hline
      71-89 & 13 & 15 \\\hline
      97 & 14 & 16 \\\hline
      131 & 16 & 17 \\\hline
      257 & 21 & 22 \\\hline
    \end{tabular}
     & 
    \begin{tabular}{|c|c|c|}
      \hline
      \multicolumn{3}{|c|}{$p=7$} \\
      \hline
      $\ell$ & Obs. & The. \\\hline\hline
      11 & 4 & 5 \\\hline
      13 & 5 & 6 \\\hline
      17 & 6 & 6 \\\hline
      19-23 & 6 & 7 \\\hline
      29-31 & 6 & 8 \\\hline
      37-61 & 8 & 10 \\\hline
      67-73 & 10 & 11 \\\hline
      79-83 & 10 & 12 \\\hline
      89-97 & 11 & 13 \\\hline
      131 & 13 & 14 \\\hline
      257 & 15 & 16 \\\hline
    \end{tabular}
     & 
    \begin{tabular}{|c|c|c|}
      \hline
      \multicolumn{3}{|c|}{$p=11$} \\
      \hline
      $\ell$ & Obs. & The. \\\hline\hline
      13 & 3 & 4 \\\hline
      17-29 & 4 & 5 \\\hline
      31 & 5 & 6 \\\hline
      37-59 & 6 & 7 \\\hline
      61 & 6 & 8 \\\hline
      67-89 & 7 & 9 \\\hline
      97 & 8 & 10 \\\hline
      131 & 9 & 11 \\\hline
      257 & 12 & 12 \\\hline
    \end{tabular}
  \end{tabular}
  \caption{$p$-adic precisions for $p=5,7,11$ and $\ell\leqslant257$.}
  \label{fig:isogenie_precision}
\end{figure}

It turns out that the precision observed in practice is near the
theoretical bound. For many values of $\ell$, a gap between the theoretical
bound and the observed bound appears, but this difference remains quite small.

\subsection{Example}
 
Let $E:y^2 = x^3 + x + 4$ be defined over $\FF_5$ and $\ell = 11$. 

We first need to compute an upper bound for the $5$-adic precision,
\begin{equation*}
  \begin{array}{ccc}
    \lploss(5, 11, 1) = 0, &  \lploss(5, 11, 2) = 1, &  \lploss(5, 11, 3) = 1, \\
    \lploss(5, 11, 4) = 2, &  \lploss(5, 11, 5) = 1\,.
  \end{array}
\end{equation*}
We find $\loss(5, 11) = 5$ and the $5$-adic precision is thus $6$. 

A $5$-adic lift of the curve is $y^2 = x^3 + x + 4 $.  With the help
of the classical $5$-th modular polynomial $\Phi_{11}$, we find that a
$11$-isogenous curve is given by $y^2 = x^3 - 7329 x - 3934 + O(5^6)$.

We can now compute the series $\bar{G}(x)$ modulo $x^{4 \ell - 1}$.
\begin{small}
\begin{multline*}
  \bar{G}(x)  =   4374 x^{42}  +  4298  x^{40}
    -  2331  x^{38}  -  4417  x^{36}
    +  3936  x^{34}  +  3505  x^{32}\\
    +  228  x^{30}  -  1041  x^{28}
    -  616  x^{26}  +  97  x^{24}
    +  236  x^{22}  +  95  x^{20}
    -  48  x^{18}  \\ -  47  x^{16}
    -  12  x^{14}  +  15  x^{12}
    +  8  x^{10}  +  x^{8}
    -  4  x^{6}  -  x^{4}
    +  1  + O(5^6) \bmod x^{43} \,.
\end{multline*}
\end{small}

A solution of the differential equation based on $\bar{G}(x)$ and $\bar{H}(z)
= \bar{a}'_6 \, z^6 \, + \, \bar{a}'_4 \, z^4 \, + \, 1$ is then given modulo
$x^{44}$ by
\begin{small}
\begin{multline*}
  \bar{S}(x)  =  - \left( 2 + O(5) \right) x^{43}  + \left( 2 + O(5) \right)
  x^{41}  - \left( 1 + O(5) \right) x^{39}  + \left( 8 + O(5^2) \right) x^{37}\\
  - \left( 1 + O(5) \right) x^{35}  + \left( O(5^2) \right) x^{33}  + \left(O(5^2)
\right) x^{31}  - \left( 10 + O(5^2) \right) x^{29}  - \left( 7 + O(5^2)
\right) x^{27}\\  - \left( 1 + O(5^2) \right) x^{25}  + \left( 192 + O(5^4)
\right) x^{23}  + \left( 125 + O(5^4) \right) x^{21}  + \left( 293 + O(5^4)
\right) x^{19} \\ + \left( 4 + O(5^4) \right) x^{17}  - \left( 161 + O(5^4)
\right) x^{15}  - \left( 611 + O(5^5) \right) x^{13}  + \left( 211 + O(5^5)
\right) x^{11}\\  - \left( 1494 + O(5^5) \right) x^{9}  + \left( 1058 + O(5^5)
\right) x^{7}  - \left( 733 + O(5^5) \right) x^{5}  + \left( O(5^6) \right)
x^{3}  + \left( 1 + O(5^6) \right) x  \,,
\end{multline*}  
\end{small}
and modulo $5$, we find
\begin{small}
\begin{multline*}
  T(x)  =  3 \, x^{21} \, + \, 2 \, x^{20} \, + \, 4 \, x^{19} \, + \, 3 \, x^{18} \, + \, 4 \, x^{17} \, + \, 3 \, x^{15} \, + \, 3 \, x^{13} \, + \, 4 \, x^{12} \, + \, 2 \, x^{11} \\
   + \, 3 \, x^{9} \, + \, 4 \, x^{8} \, + \, 4 \, x^{7} \, + \, 4 \, x^{6} \,
   + \, x^{5} \, + \, x^{4} \, + \, 3 \, x^{3} \, + \, 2 \, x^{2} \, + \, 1
   \bmod x^{22} \,.
\end{multline*}
\end{small}
We have $R(x)=1/T(x)^2 \bmod x^{2 \ell}$, that is
\begin{small}
\begin{multline*}
    R(x)  =  2 \, x^{20} \, + \, 2 \, x^{19} \, + \, 3 \, x^{18} \, + \, x^{16} \, + \, 2 \, x^{15} \, + \, 3 \, x^{14} \, + \, x^{13} \, + \, 3 \, x^{12} \, + \, 2 \, x^{11} \\
    \, + \, 2 \, x^{10} + \, 2 \, x^{8} \, + \, 3 \, x^{7} \, + \, 4 \, x^{6} \, + \, 4 \, x^{5} \, + \, 4 \, x^{3} \, + \, x^{2} \, + \, 1 \bmod x^{22} \,.
\end{multline*}
\end{small}
The rebuilding of the rational fraction corresponding to $R$ gives
\begin{small}
  \begin{displaymath}
    R(x) \, = \, \frac{3 \, x^{11} \, + \, x^{9} \, + \, x^{8} \, + \, x^{7} \, + \, x^{6} \, + \, 3 \, x^{5} \, + \, 2 \, x^{4} \, + \, 3 \, x^{3} \, + \, 2 \, x^{2} \, + \, 2 \, x \, + \, 1}{x^{10} \, + \, x^{9} \, + \, x^{8} \, + \, x^{7} \, + \, 3 \, x^{6} \, + \, 3 \, x^{5} \, + \, 3 \, x^{4} \, + \, 2 \, x^{3} \, + \, x^{2} \, + \, 2 \, x \, + \, 1} \bmod x^{22}.
  \end{displaymath}
\end{small}
One reverses the order of the coefficients of the denominator to obtain
\begin{small}
\begin{displaymath}
  D(x) \, = \, x^{10} \, + \, 2 \, x^9 \, + \, x^8 \, + \, 2 \, x^7 \, + \, 3 \, x^6 \, + \, 3 \, x^5 \, + \, 3 \, x^4 \, + \, x^3 \, + \, x^2 \, + \, x \, + \, 1.
\end{displaymath}
\end{small}
The $\ell$-th Elkies polynomial is then
\begin{small}
  \begin{displaymath}
    \sqrt{D(x)} \, = \, x^5 \, + \, x^4 \, + \, x^2 \, + \, 3 \, x \, + \, 1.
  \end{displaymath}
\end{small}

\iffalse
\bibliographystyle{plain}
\bibliography{padiso2}
\else

\fi

\appendix
\section{Proof of Proposition~\ref{prop:eqdiff}}
\label{sec:proof}

Let $d$ be a non-zero even integer, we assume that we know a
solution of the differential equation modulo $x^{d+1}$. We thus have
\begin{equation}
  \label{isogenie_Sd}
  {S_d'}^2 = G \cdot (H \circ S_d) \bmod x^d \, , \, S_d(0) = \alpha \, , \, S_d'(0) = \beta.
\end{equation}
Let  $S_{2d} = S_d + A_{2d}$ be a solution modulo $x^{2d+1}$, with $x^{d+1}$
dividing $A_{2d}$, therefore 
\begin{math}
  \left( S_d' + A_{2d}' \right)^2 = G \cdot \big( H \circ (S_d + A_{2d}) \big) \bmod x^{2d}\,.
\end{math}
This yields a linear differential equation in $A_{2d}$.
\begin{displaymath}
  2 \, S_d' \cdot A_{2d}' - G \cdot (H' \circ S_d) \cdot A_{2d} = G \cdot (H
  \circ S_d) - {S_d'}^2 \bmod x^{2d}\,. 
\end{displaymath}
With initial condition $A_{2d}(0) = 0$, a solution of this equation is
\begin{equation}
  \label{isogenie_A2d}
  A_{2d} = \frac{1}{J_{2d}} \, \int \frac{\big( G \cdot (H \circ S_d) - {S_d'}^2 \big) \cdot J_{2d}}{2 \, S_d'} \; dx \bmod x^{2d+1} \, ,
\end{equation}
\begin{displaymath}
  \text{where} \quad J_{2d} = \exp \left( - \int \frac{G \cdot (H' \circ
      S_d)}{2 \, S_d'} \; dx \right) \bmod x^{2d+1} \,.
\end{displaymath}

From Eq.~(\ref{isogenie_Sd}), we know that $\big( G \cdot (H \circ S_d) -
{S_d'}^2 \big)$ is divisible by $x^d$. Moreover, $S_d'$ has a non-zero
constant coefficient. A factor $x^d$ appears then in the integral and it's
enough to compute $J_{2d}$ modulo $x^d$.
The inverse of $J_{2d}$ is multiplied by the integral, it will thus be multiplied by
$x^{d+1}$, and it's enough to evaluate this inverse modulo $x^d$.
The inverse of $S_d'$ is needed in the computations of $A_{2d}$ and $J_{2d}$.
In $A_{2d}$, this inverse is multiplied by $x^d$ and we then compute a
primitive.  In $J_{2d}$, we compute only modulo $x^d$. In both cases, the
inverse of $S'_d$ modulo $x^d$ is enough. This inverse is provided by
Eq.~(\ref{isogenie_Sd}):
\begin{displaymath}
  \frac{1}{S'_d} = \frac{S'_d}{G \cdot (H \circ S_d)} \bmod x^d\,.
\end{displaymath}
We plug this expression in the computation of $J_{2d}$ modulo $x^d$, we find
\begin{eqnarray*}
  \int \frac{G \cdot (H' \circ S_d)}{2 \, S'_d} \, dx & = & \int \frac{S'_d \cdot (H' \circ S_d)}{2 \, (H \circ S_d)} \, dx \bmod x^d \\
  & = & \frac{\log(H \circ S_d)}{2} \bmod x^d.
\end{eqnarray*}
We then find the following nice formulas for $J_{2d}$ and $1 / J_{2d}$ modulo $x^d$,
\begin{displaymath}
  J_{2d} = \frac{1}{\sqrt{H \circ S_d}} \bmod x^d \, , \quad \frac{1}{J_{2d}} = \sqrt{H \circ S_d} \bmod x^d \,.
\end{displaymath}

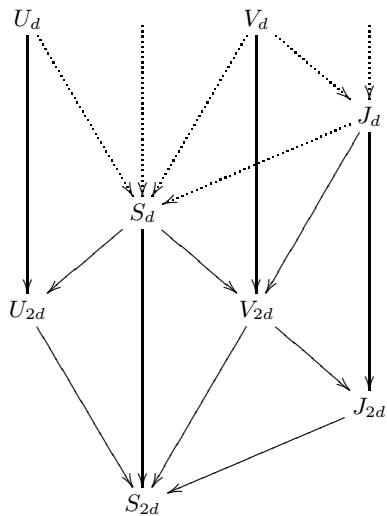
\begin{figure}[htb]
  \centering
    \mbox{\small
      \xymatrix{
        U_d \ar[ddd] \ar@{.>}[ddr] & \ar@{.>}[dd] & V_d \ar@{.>}[ddl] \ar[ddd] \ar@{.>}[dr] & \ar@{.>}[d] \\
        & & & J_d \ar@{.>}[dll] \ar[ddl] \ar[ddd] \\
        & S_d \ar[dl] \ar[ddd] \ar[dr] \\
        U_{2d} \ar[ddr] & & V_{2d} \ar[ddl] \ar[dr] \\
        & & & J_{2d} \ar[dll] \\
        & S_{2d} \\
      }
    }
    \caption{Computation of $U_{2d}$, $V_{2d}$, $J_{2d}$ and $S_{2d}$}
    \label{fig:eqdif}
\end{figure}
These formulas allow to efficiently compute $S_{2d}$ from $S_d$ and other
known quantities.
\begin{itemize}
\item From the inverse of $S'_{d/2}$ modulo $x^{d/2}$, denoted by $U_d$, we use a
  classical Newton iteration to compute $U_{2d}$. Since $S_d = S_{d/2}
  \bmod x^{d/2+1}$, we have $U_{2d} = U_d \bmod x^{d/2}$ and we compute the
  coefficients of $U_{2d}$ thanks to
  \begin{displaymath}
    U_{2d} = U_d \cdot \big( 2 - S'_d \cdot U_d \big) \bmod x^d\,.
  \end{displaymath}
\item From $\sqrt{H \circ S_{d/2}}$ modulo $x^{d/2}$, denoted by $V_d$, and the inverse of $V_d$ modulo $x^{d/2}$, denoted by $J_d$, we compute $V_{2d}$ and $J_{2d}$ as follows. Getting $V_{2d}$ consists in computing a solution of $v^2 - (H \circ S_d)(x) = 0$. We use
  \begin{displaymath}
    V_{2d} = \frac{1}{2} \left( V_d +  \frac{H \circ S_d}{V_d} \right) \bmod x^d\,.
  \end{displaymath}
  $J_d$ and $V_d$ are by definition inverses of each other modulo $x^{d/2}$. We obtain the inverse $W_{2d}$ of $V_d$ modulo $x^d$ thanks to Newton formulas too,
  \begin{displaymath}
    W_{2d} = J_d \cdot \big( 2 - V_d \cdot J_d \big) \bmod x^d\,.
  \end{displaymath}
  If we now plug this value in the $V_{2d}$ formula, we finally find
  \begin{displaymath}
    2\, V_{2d} = V_d + J_d \cdot (H \circ S_d) \cdot (2 - V_d \cdot J_d) \bmod x^d.
  \end{displaymath}
  Another use of Newton's inversion formula yields $J_{2d}$,
  \begin{displaymath}
    J_{2d} = J_d \cdot \big( 2 - J_d \cdot V_{2d} \big) \bmod x^d\,.
  \end{displaymath}
\end{itemize}

Thanks to all these equations, we can compute $(U_{2d}, V_{2d}, J_{2d})$ from
$(U_d, V_d, J_d, S_d)$. The quantity $S_{2d}$ is then obtained from Eq.~(\ref{isogenie_A2d}),
\begin{displaymath}
  S_{2d} = S_d + \frac{V_{2d}}{2} \, \int U_{2d} \cdot J_{2d} \cdot \big( G
  \cdot (H \circ S_d) - S_d'^2 \big) \, dx \bmod x^{2d+1}.
\end{displaymath}
We illustrate the corresponding computations in Fig.~\ref{fig:eqdif}. \medskip

It remains to obtain initial values, for $d = 2$. Let $\gamma$ be defined by
$S_2(x) = \alpha + \beta \, x + \gamma \, x^2 \bmod x^3$. The series $S_2$ is
solution of the differential equation modulo $x^2$ and thus $\beta^2 + 4 \,
\beta \, \gamma \, x = G(x) \, H(\alpha + \beta x) \bmod x^2$.
Once derivated, and evaluated at $x=0$, we obtain $\gamma$, and thus the value
of $S_2$,
\begin{displaymath}
  S_2(x) = \alpha + \beta \, x + \left( \frac{G'(0)}{4 \beta} + \frac{\beta^2
      \, H'(\alpha)}{4} \right) \, x^2 \bmod x^3\,.
\end{displaymath}
We deduce as well
\begin{equation*}
  U_2(x) = \frac{1}{\beta} \bmod x \, , \quad V_2(x) = 1 \bmod x \quad \text{
    and } \quad J_2(x) = 1 \bmod x.
\end{equation*}

\medskip

\end{document}